\documentclass[12pt,twoside,reqno]{amsart}
\usepackage[colorlinks=false,citecolor=blue]{hyperref}
\usepackage{mathptmx, amsmath, amssymb, amsfonts, amsthm, mathptmx, enumerate, color,mathrsfs}
\usepackage{graphicx}

\usepackage{epstopdf}
\newtheorem{theorem}{Theorem}[section]
\newtheorem{lemma}{Lemma}[section]

\newtheorem{definition}{Definition}[section]
\newtheorem{corollary}{Corollary}[section]
\newtheorem{example}{Example}[section]

\numberwithin{equation}{section}

\begin{document}
	
\title{Strongly generalized derivations on $C^{\ast}$-algebras}
\author{Amin Hosseini}

\begin{abstract}
Let $\mathcal{A}$ and $\mathcal{B}$ be two algebras, let $\mathcal{M}$ be a $\mathcal{B}$-bimodule and let $n$ be a positive integer. A linear mapping $D_n:\mathcal{A} \rightarrow \mathcal{M}$ is called a strongly generalized derivation of order $n$, if there exist the families $\{E_k:\mathcal{A} \rightarrow \mathcal{M}\}_{k = 1}^{n}$, $\{H_k:\mathcal{A} \rightarrow \mathcal{M}\}_{k = 1}^{n}$, $\{F_k:\mathcal{A} \rightarrow \mathcal{B}\}_{k = 1}^{n}$ and $\{G_k:\mathcal{A} \rightarrow \mathcal{B}\}_{k = 1}^{n}$ of mappings which satisfy $$D_n(ab) = \sum_{k = 1}^{n}\left[E_k(a) F_k(b) + G_k(a)H_k(b)\right]$$ for all $a, b \in \mathcal{A}$. In this paper, we prove that every strongly generalized derivation of order one from a $C^{\ast}$-algebra into a Banach bimodule is automatically continuous under certain conditions. The main theorem of this paper extends some celebrated results in this regard.
\end{abstract}

\maketitle

\pagestyle{myheadings}
\markboth{\centerline {}}
{\centerline {}}
\bigskip
\bigskip

\section{Introduction and Preliminaries}
Let $\mathcal{A}$ and $\mathcal{B}$ be two algebras, let $\mathcal{M}$ be a $\mathcal{B}$-bimodule and let $n$ be a positive integer. A linear mapping $D_n:\mathcal{A} \rightarrow \mathcal{M}$ is called a strongly generalized derivation of order $n$, if there exist the families $\{E_k:\mathcal{A} \rightarrow \mathcal{M}\}_{k = 1}^{n}$, $\{H_k:\mathcal{A} \rightarrow \mathcal{M}\}_{k = 1}^{n}$, $\{F_k:\mathcal{A} \rightarrow \mathcal{B}\}_{k = 1}^{n}$ and $\{G_k:\mathcal{A} \rightarrow \mathcal{B}\}_{k = 1}^{n}$ of mappings which satisfy $$D_n(ab) = \sum_{k = 1}^{n}\left[E_k(a) F_k(b) + G_k(a)H_k(b)\right]$$ for all $a, b \in \mathcal{A}$. Clearly, for $n = 1$, we have
$$D_1(ab) = E(a) F(b) + G(a)H(b)$$ for all $a, b \in \mathcal{A}$, where $E, H:\mathcal{A} \rightarrow \mathcal{M}$ and $F, G:\mathcal{A} \rightarrow \mathcal{B}$ are arbitrary mappings. If $D_1:\mathcal{A} \rightarrow \mathcal{M}$ is a strongly generalized derivation of order one associated with the mappings $E, H: \mathcal{A} \rightarrow \mathcal{M}$ and $F, G:\mathcal{A} \rightarrow \mathcal{B}$, then we say that $D_1$ is an $(E, F, G, H)$-derivation. Also, if $D_n$ is a strongly generalized derivation of order $n$ associated with the families $\{E_k:\mathcal{A} \rightarrow \mathcal{M}\}_{k = 1}^{n}$, $\{H_k:\mathcal{A} \rightarrow \mathcal{M}\}_{k = 1}^{n}$, $\{F_k:\mathcal{A} \rightarrow \mathcal{B}\}_{k = 1}^{n}$ and $\{G_k:\mathcal{A} \rightarrow \mathcal{B}\}_{k = 1}^{n}$ of mappings, we say that $D_n$ is an $\big(\{E_k\}_{k = 1}^{n}, \{F_k\}_{k = 1}^{n}, \{G_k\}_{k = 1}^{n}, \{H_k\}_{k = 1}^{n}  \big)$-derivation. As can be seen, if $D_1$ is a strongly generalized derivation of order one, then it covers the notion of a derivation (if $D_1 = E = H$ and $F = G = I$), the notion of a generalized $(\sigma, \tau)$-derivation associated with a mapping $d$ (if $D_1 = E$, $F = \sigma$, $G = \tau$ and $H = d$), the notion of a left $\sigma$-centralizer (if $D_1 = E$, $F = \sigma$ and $G$ or $H$ is zero), the notion of a right $\tau$-centralizer (if $E$ or $F$ is zero, $G = \tau$ and $H = D_1$), the notion of a generalized derivation associated with a mapping $d$ (if $D_1 = E$, $F = G = I$ and $H = d$), the notion of a homomorphism (if $D_1 = E = F$ and $G = 0$ or $H = 0$), and the notion of a ternary derivation (if $F = G = I$). Also, if $D_2$ is a strongly generalized derivation of order two, we have
\begin{align*}
D_2(ab) = E_1(a) F_1(b) + G_1(a)H_1(b) + E_2(a) F_2(b) + G_2(a)H_2(b)
\end{align*}
for all $a, b \in \mathcal{A}$, where $E_i, H_i:\mathcal{A} \rightarrow \mathcal{M}$ and $F_i, G_i:\mathcal{A} \rightarrow \mathcal{B}$
are mappings for any $i \in \{1,2\}$.  For example, every $(\delta, \varepsilon)$-double derivation is a strongly generalized derivation of order two. For more material about $(\delta, \varepsilon)$-double derivations, see, e.g. \cite{M1}. Now we present an example of a strongly generalized derivation of order $n$. Let $\mathcal{A}$ and $\mathcal{B}$ be two algebras. A sequence $\{f_n\}$ of linear mappings from $\mathcal{A}$ into $\mathcal{B}$ is called a ternary higher derivation associated with the sequences $\{g_n\}$ and $\{h_n\}$ of mappings from $\mathcal{A}$ into $\mathcal{B}$, which is denoted by $(f_n, g_n, h_n)$, if $$f_n(ab) = \sum_{k = 0}^{n}g_{n - k}(a) h_k(b)$$ holds for all $a, b \in \mathcal{A}$ and all nonnegative integers $n$. Let $n$ be a positive integer and let $(f_n, g_n, h_n)$ be a ternary higher derivation. Then every $f_n$ is a strongly generalized derivation of order  $m$ in which
\[m = \left\{{\begin{array}{*{20}{c}}
\begin{array}{l}
\frac{n + 2}{2} \\
\frac{n + 1}{2}
\end{array}&
\begin{array}{l}
n \ is \ even,\\
n \ is \ odd
\end{array}
\end{array}} \right.\]
For more details about the structure of higher derivations, see, e.g. \cite{H7, H8} and the references therein.
It is interesting to note that the applications of generalized types of derivations, such as generalized derivations and $(\sigma, \tau)$-derivations to important physical topics have been recently studied. See, for example, \cite{Hel} for the application of generalized derivations in general relativity, and \cite{E, H5*} for the application of $(\sigma, \tau)$-derivations in theoretical physics. Therefore, it is possible that the notion of strongly generalized derivation of order $n$ be considered by physicists in the future and used in the study of physical topics. So, it seems interesting to investigate the details of these mappings. The main purpose of this paper is to investigate the automatic continuity of strongly generalized derivations of order one on $C^{\ast}$-algebras.
\\
Derivations and their various properties are significant subjects in the study of Banach algebras and $C^{\ast}$-algebras. Let $\mathcal{A}$ be a Banach or $C^{\ast}$-algebra and let $\mathcal{M}$ be a Banach $\mathcal{A}$-bimodule. One of the most important problems related to these mappings is the question that under what conditions is a derivation $d:\mathcal{A} \rightarrow \mathcal{M}$ continuous? This question lies in the theory of automatic continuity which is an important subject in operator theory and mathematical analysis and also has attracted the attention of researchers during the last few decades. In this theory, we are looking for conditions which guarantee that a linear mapping between two Banach algebras (or two Banach spaces, in general) is necessarily continuous. It is worth to note that there is an extensive literature on this topic, and we try to give a brief background in this regard. In 1958, Kaplansky \cite{K} conjectured that every derivation on a $C^{\ast}$-algebra is continuous. Two years later, Sakai \cite{S} answered this conjecture. Indeed, he proved that every derivation on a $C^{\ast}$-algebra is automatically continuous and later in 1972, Ringrose \cite{R}, by using the pioneering work of Bade and Curtis \cite{B1} concerning the automatic continuity of a module homomorphism between bimodules over $C(K)$-spaces, showed that every derivation from a $C^{\ast}$-algebra $\mathcal{A}$ into a Banach $\mathcal{A}$-bimodule is automatically continuous. Another celebrated theorem in this regard, proven by Johnson and Sinclair \cite{Ja}, states that every derivation on a semisimple Banach algebra is continuous. In addition, in an interesting article, Peralta and Russo \cite{P} investigated automatic continuity of derivations on $C^{\ast}$-algebras and $JB^{\ast}$-triples. We further know from \cite{EP18} that every generalized derivation on a von Neumann algebra and every linear mapping on a von Neumann algebra which is a derivation or a triplet derivation at zero is automatically continuous. Recently, the present author, in collaboration with Peralta and Su \cite{H0}, presented important and interesting results about the continuity of generalized derivations and ternary derivations on $C^{\ast}$-algebras. In addition, the author of this article has studied the continuity of $(\delta, \varepsilon$)-double derivations, $(\sigma, \tau)$-derivations and $\phi$-derivations on $C^{\ast}$-algebras and Banach algebras, see \cite{H3, H4, H5, H6}. Moreover, Hou and Ming \cite{Hu} proved that if $\mathcal{X}$ is simple and $\sigma, \tau$ are surjective and continuous mappings on $B(\mathcal{X})$, then every $(\sigma, \tau)$-derivation on $B(\mathcal{X})$ is continuous, where $B(\mathcal{X})$ denotes the algebra of all bounded linear mappings from $\mathcal{X}$ into itself. We refer the reader to \cite{D, D1, Ru, S1, V1, V2, V3} for a deep and extensive study on this subject. Now we turn to the main theorem of this article.\\
Let $\mathcal{A}$ be a $C^{\ast}$-algebra, let $\mathcal{B}$ be a Banach algebra, let $\mathcal{M}$ be a Banach $\mathcal{B}$-bimodule and let $D_1: \mathcal{A} \rightarrow \mathcal{M}$ be an $(E, F, G, H)$-derivation such that $E, H: \mathcal{A} \rightarrow \mathcal{M}$ are linear and $F, G:\mathcal{A} \rightarrow \mathcal{B}$ are continuous mappings at zero. Moreover, assume that either $(G(ab) - G(a)G(b))H(c) = 0$ or $E(c) (F(ab) - F(a)F(b)) = 0$ for all $a, b, c \in \mathcal{A}$. Then $D_1$ is continuous.\\
Moreover, some consequences of the abovementioned result are presented. In fact, by considering the notion of an $(E, F, G, H)$-derivation, we obtain the above-mentioned result for derivations, generalized $(\sigma, \tau)$-derivations, left (resp. right) centralizers, generalized derivations and ternary derivations.

\section{Results and Proofs}
First of all, we give some examples of $(E, F, G, H)$-derivations.


\begin{example}Let $\mathcal{A}$ be an algebra, and let
\begin{align*}
\mathfrak{A} = \Bigg\{\left [\begin{array}{ccc}
a & b & c\\
0 & 0 & 0\\
0 & 0 & e
\end{array}\right ] \ : \ a, b, c, e \in \mathcal{A}\Bigg\}
\end{align*}
Clearly, $\mathfrak{A}$ is an algebra under the usual matrix operations. Define the mappings $D_1, E, F, G, H:\mathfrak{A} \rightarrow \mathfrak{A}$ by $$D\Bigg(\left [\begin{array}{ccc}
a & b & c\\
0 & 0 & 0\\
0 & 0 & e
\end{array}\right ]\Bigg) = \left [\begin{array}{ccc}
0 & a & b\\
0 & 0 & 0\\
0 & 0 & 0
\end{array}\right ],$$

$$ E\Bigg(\left [\begin{array}{ccc}
a & b & c\\
0 & 0 & 0\\
0 & 0 & e
\end{array}\right ]\Bigg) = \left [\begin{array}{ccc}
a & b & 0\\
0 & 0 & 0\\
0 & 0 & e
\end{array}\right ],$$

$$ H\Bigg(\left [\begin{array}{ccc}
a & b & c\\
0 & 0 & 0\\
0 & 0 & e
\end{array}\right ]\Bigg) = \left [\begin{array}{ccc}
-a & 0 & b\\
0 & 0 & 0\\
0 & 0 & e
\end{array}\right ],$$

$$ F\Bigg(\left [\begin{array}{ccc}
a & b & c\\
0 & 0 & 0\\
0 & 0 & e
\end{array}\right ]\Bigg) = \left [\begin{array}{ccc}
a & a & 0\\
0 & 0 & 0\\
0 & 0 & e
\end{array}\right ],$$

$$ G\Bigg(\left [\begin{array}{ccc}
a & b & c\\
0 & 0 & 0\\
0 & 0 & e
\end{array}\right ]\Bigg) = \left [\begin{array}{ccc}
a & 0 & 0\\
0 & 0 & 0\\
0 & 0 & -e
\end{array}\right ].$$

A simple calculation shows that \begin{align*}
D_1(AB) = E(A)F(B) + G(A) H(B), \ \ \ \ \ \ \ \ \ \ \ A, B \in \mathfrak{A},
\end{align*}
which means that $D_1$ is an $(E, F, G, H)$-derivation on $\mathfrak{A}$.
\end{example}

\begin{example} Let $\mathcal{A}$ and $\mathcal{B}$ be two algebras. It is easy to see that $\mathfrak{A} = \mathcal{A} \times \mathcal{B}$ is an algebra by the following product:
$$(a_1, b_1) \bullet (a_2, b_2) = (a_1 a_2, b_1 b_2)$$
for all $a_1, a_2 \in \mathcal{A}$ and $b_1, b_2 \in \mathcal{B}$. Let $f:\mathcal{A} \rightarrow \mathcal{A}$ and $g:\mathcal{B} \rightarrow \mathcal{B}$ are two mappings (linear or nonlinear). Define the mappings $D_1, E, F, G, H:\mathfrak{A} \rightarrow \mathfrak{A}$ by
\begin{align*}
&D_1((a,b)) = (a, 0), \\ & E((a,b)) = (a, g(b)), \\ & H((a,b)) = (0, -\lambda_1 a - \lambda_2 b), \\ & F((a,b)) = (a,\lambda_1 a + \lambda_2 b), \\ & G((a,b)) = (f(a), g(b)),
\end{align*}
where $\lambda_1$ and $\lambda_2$ are two complex numbers. A routine calculation shows that $D_1$ is an $(E, F, G, H)$-derivation on $\mathfrak{A}$.
\end{example}

We state the following auxiliary Lemmas which will be used extensively to prove the main theorem of this paper.
\begin{lemma} \label{1} Let $\mathcal{\mathfrak{I}}$ be a closed ideal in a $C^{\ast}$-algebra $\mathcal{A}$ and let $\{a_1, a_2, a_3, ...\}$ be a subset of $\mathcal{\mathfrak{I}}$ such that $\sum_{n = 1}^{\infty}\|a_n\|^{2} \leq 1$. Then there exist elements $b, c_1, c_2, ...$ of $\mathcal{\mathfrak{I}}$ such that $b \geq 0$, $\|c_n\| \leq 1$, and $a_n = b c_n$ for any $n \in \mathbb{N}$.
\end{lemma}
\begin{proof} See \cite[Exercise 4.6.40]{K-R}.
\end{proof}

\begin{lemma} \label{2} Suppose that $\mathcal{A}$ is an infinite-dimensional $C^{\ast}$-algebra. Then there is an infinite sequence $\{a_1, a_2, a_3, ...\}$ of nonzero elements of $\mathcal{A}^{+}$ such that $a_j a_k = 0$ when $j \neq k$.
\end{lemma}
\begin{proof} See \cite[Exercise 4.6.13]{K-R}.
\end{proof}

\begin{lemma} \label{3} Let $\mathcal{A}$ and $\mathcal{B}$ be two $C^{\ast}$-algebras, and let $\varphi$ be a $\ast$-homomorphism from $\mathcal{A}$ onto $\mathcal{B}$. Let $\{b_1, b_2, b_3, ...\}$ be a sequence of elements of $\mathcal{B}^{+}$ such that $b_j b_k = 0$ when $ j \neq k$. Then there is a sequence $\{a_1, a_2, a_3, ...\}$ of elements of $\mathcal{A}^{+}$ such that $a_j a_k =0$ when $j \neq k$, and $\varphi(a_j) = b_j$ for any $j \in \mathbb{N}$.
\end{lemma}
\begin{proof} See \cite[exercise 4.6.20]{K-R}.
\end{proof}

We are now ready to prove the main result of the this paper, which is inspired by the proof of the related results on the ordinary derivations and $(\sigma, \tau)$-derivations in \cite{R} and \cite{Hu}, respectively.

\begin{theorem} \label{4} Let $\mathcal{A}$ be a $C^{\ast}$-algebra, let $\mathcal{B}$ be a Banach algebra, let $\mathcal{M}$ be a Banach $\mathcal{B}$-bimodule and let $D_1: \mathcal{A} \rightarrow \mathcal{M}$ be an $(E, F, G, H)$-derivation such that $E, H: \mathcal{A} \rightarrow \mathcal{M}$ are linear and $F, G:\mathcal{A} \rightarrow \mathcal{B}$ are continuous mappings at zero. Moreover, assume that either $(G(ab) - G(a)G(b))H(c) = 0$ or $E(c) (F(ab) - F(a)F(b)) = 0$ for all $a, b, c \in \mathcal{A}$. Then $D_1$ is continuous.
\end{theorem}

\begin{proof} Suppose that $D_1:\mathcal{A} \rightarrow \mathcal{M}$ is an $(E, F, G, H)$-derivation such that $F, G:\mathcal{A} \rightarrow \mathcal{A}$ are continuous mappings and $(G(ab) - G(a)G(b))H(c) = 0$ for all $a, b, c \in \mathcal{A}$. Using a five-step proof, we show that $D_1$ is continuous. For each $a \in \mathcal{A}$, we consider the mappings $\gamma_a : \mathcal{A} \rightarrow \mathcal{M}$, $\gamma_a (t) = D_1(at)$ and
$\psi_a : \mathcal{A} \rightarrow \mathcal{M}$, $\psi_a (t) = G(a)H(t)$. Before entering the first step of our proof, note that the mapping $\Lambda_a : \mathcal{A} \rightarrow \mathcal{M}$ defined by $\Lambda_a(t) = E(a)F(t)$, where $a$ is an arbitrary element of $\mathcal{A}$, is linear. To see this, simply use the linearity of $D_1$ and $H$. We leave the details to the interested reader. Similarly, the linearity of $D_1$ and $E$ implies the linearity of the mapping $\Omega_a:\mathcal{A} \rightarrow \mathcal{M}$ defined by $\Omega_a(t) = G(t) H(a)$, where $a$ is an arbitrary element of $\mathcal{A}$. Since the mappings $F$ and $G$ are continuous at zero, the linear mappings $\Lambda_a$ and $\Omega_a$ are continuous.
Let $\mathcal{I} = \{a \in \mathcal{A} \ : \ \gamma_a \ is \ continuous\}$.\\
\\
\textbf{Step 1:} $\mathcal{I} = \{a \in \mathcal{A} \ : \ \psi_a \ is \ continuous\}$.\\
Set $\mathcal{V} = \{a \in \mathcal{A} \ | \ \psi_a \ is \ continuous\}$. Our task is to show that $\mathcal{V} = \mathcal{I}$. Let $a$ be an element of $\mathcal{I}$. It is clear that the mapping $t \mapsto D_1(at) - E(a)F(t) = G(a) H(t)$ is continuous, which means that $a \in \mathcal{V}$ and thus, $\mathcal{I} \subseteq \mathcal{V}$. Now, we prove that $\mathcal{V} \subseteq \mathcal{I}$. Let $a$ be an element of $\mathcal{V}$. So, the mapping $t \mapsto G(a)H(t)$ is continuous. Therefore, the mapping $t \mapsto E(a)F(t) + G(a)H(t) = D_1(at)$ is continuous. This yields that $a \in \mathcal{I}$ and it means that $\mathcal{V} \subseteq \mathcal{I}$. Consequently, $\mathcal{I} = \{a \in \mathcal{A} \ : \ \psi_a \ is \ continuous\}$.
\\
\\
\textbf{Step 2}: $\mathcal{I}$ is a closed two sided-ideal of $\mathcal{A}$.\\ \\
First, we show that $\mathcal{I}$ is a two sided-ideal of $\mathcal{A}$. Note that for every element $b \in \mathcal{A}$, the linear mapping $\theta:\mathcal{A} \rightarrow \mathcal{A}$ defined by $\theta(t) = bt$ is continuous. Assume that $a$ and $b$ are two arbitrary elements of $\mathcal{I}$ and $\mathcal{A}$, respectively. It is evident that the mapping $\gamma_a o \theta:\mathcal{A} \rightarrow \mathcal{M}$ defined by $\gamma_a o \theta(t) = D_1(ab t)$ is continuous, and so $ab \in \mathcal{I}$. It means that $\mathcal{I}$ is a right ideal of $\mathcal{A}$. Since $a \in \mathcal{I} = \mathcal{V}$, the mapping $t \mapsto G(b)G(a) H(t)$ is continuous. Obviously, the mapping $t \mapsto E(ba)F(t)$ is continuous and since we are assuming that $(G(ba) - G(b)G(a))H(t) = 0$ for all $a, b, t \in \mathcal{A}$, the mapping $t \mapsto E(ba)F(t) + G(ba) H(t) = D_1(bat)$ is continuous. This yields that $ba \in \mathcal{I}$ and consequently, $\mathcal{I}$ is a left ideal of $\mathcal{A}$. Hence, $\mathcal{I}$ is a bi-ideal of $\mathcal{A}$. In the following, we show that $\mathcal{I}$ is closed. Let $a \in \overline{\mathcal{I}}$. Then there exists a sequence $\{a_n\} \subseteq \mathcal{I}$ such that $\lim_{n \rightarrow \infty}a_n = a$.
It is enough to show that the mapping $\psi_a:\mathcal{A} \rightarrow \mathcal{M}$ is continuous, i.e. $a \in \mathcal{V}$. Since $\{a_n\}$ is a sequence in $\mathcal{V}$, the linear mapping $\psi_{a_{n}}:\mathcal{A} \rightarrow \mathcal{M}$ is continuous for every $n \in \mathbb{N}$. We have $\lim_{n \rightarrow \infty} \psi_{a_n}(t) = \psi_a(t)$. By the principle of uniform boundedness, $\psi_a$ is norm continuous and so $a \in \mathcal{V} = \mathcal{I}$. Therefore, $\mathcal{I}$ is a closed two sided-ideal of $\mathcal{A}$.
\\
\\
\textbf{Step 3}: $D_1 |_{\mathcal{I}}$ is continuous. \\ \\
Suppose that $D_1|_{\mathcal{I}}$ is an unbounded linear mapping. It means that $\|D_1 |_{\mathcal{I}}\| = sup\{\|D_1(a_n)\| \ : \ \|a_n\| \leq 1, \ a_n \in \mathcal{I}\} = \infty$. Then, we can choose a sequence $\{a_n\}$ in $\mathcal{I}$ such that $\|D_1(a_n)\| \rightarrow \infty$, $\sum_{n = 1}^{\infty}\|a_n\|^{2} \leq 1$. Now we define $b = (\sum_{n = 1}^{\infty}a_n a_n^{\ast})^{\frac{1}{4}}$, and since $\mathcal{I}$ is a closed bi-ideal of $\mathcal{A}$, $b$ is a positive element of $\mathcal{I}$, i.e. $b \in \mathcal{I}^{+}$. We have $\|b\|^{4} = \|b^4\| = \|\sum_{n = 1}^{\infty}a_n a_n^{\ast}\| \leq \sum_{n = 1}^{\infty}\|a_n a_n^{\ast}\| = \sum_{n = 1}^{\infty}\|a_n\|^{2} \leq 1$, which implies that $\|b\| \leq 1$. It follows from Lemma \ref{1} that for every $n \in \mathbb{N}$ there exists an element $c_n \in \mathcal{I}$ such that $\|c_n\| \leq 1$, $a_n = b c_n$. Note that $\|D_1(b c_n)\| = \|D_1(a_n)\| \rightarrow \infty$. Therefore, we have $\infty = sup\{\|D_1(b c_n)\| \ : \ \|c_n\| \leq 1\} \leq sup \{\|D_1(bt)\| \ : \ \|t\| \leq 1\}$, and consequently, the mapping $\gamma_{b}:\mathcal{A} \rightarrow \mathcal{M}$ defined by $\gamma_{b}(t) = D_1(bt)$ is unbounded. But this is a contradiction of the fact that $b \in \mathcal{I}$. Hence, the restriction $D_1 |_{\mathcal{I}}$ is continuous.
\\
\\\textbf{Step 4}: $\frac{\mathcal{A}}{\mathcal{I}}$ is finite-dimensional.\\ \\
To obtain a contradiction, assume that $\frac{\mathcal{A}}{\mathcal{I}}$ is an infinite-dimensional $C^{\ast}$-algebra. It follows from Lemma \ref{2} that there exists an infinite sequence $\{b_1, b_2, b_3, ...\}$ of non-zero, positive elements in $\frac{\mathcal{A}}{\mathcal{I}}$ such that $b_j b_k = 0$ where $j \neq k$. Since $\|b_j^{2}\| = \|b_j\|^{2} > 0$, we have $b_j^{2} \neq 0$. We know that the natural mapping $\pi:\mathcal{A} \rightarrow \frac{\mathcal{A}}{\mathcal{I}}$ is a $\ast$-homomorphism from the $C^{\ast}$-algebra $\mathcal{A}$ onto the $C^{\ast}$-algebra $\frac{\mathcal{A}}{\mathcal{I}}$. According to Lemma \ref{3}, there exists a sequence $\{s_1, s_2, s_3, ...\}$ of elements of $\mathcal{A}^{+}$ such that $\pi(s_j) = b_j$, $s_j s_k = 0$, where $j \neq k$. If we now replace $s_j$ by an appropriate scalar multiple, we may suppose also that $\|s_j\| \leq 1$. It follows from $\pi(s_j^{2}) = b_j^{2} \neq 0$ that $s_j^{2} \notin \mathcal{I}$. This fact along with the definition of $\mathcal{I}$ imply that the mapping $\eta_{s_j^{2}}:\mathcal{A} \rightarrow \mathcal{M}$ defined by $t \mapsto D_1(s_j^2 t)$ is unbounded. Let $\{t_j\}$ be a sequence of $\mathcal{A}$ such that $\|t_j\| \leq 2^{-j}$. Since $\sum\|s_j t_j\| \leq \sum\|s_j\| \|t_j\| \leq \sum2^{-j} < \infty$, the series $\sum s_j t_j$ converges to an  element $c$ of $\mathcal{A}$, i.e. $\sum s_j t_j = c$.  Hence, $\|c\| = \|\sum s_j t_j\| \leq \sum \|s_j t_j\| \leq \sum 2^{-j} \leq 1$. Note that $s_j c = s_j(\sum s_j t_j) = s_j s_1 t_1 + s_j s_2 t_2 + ... + s_j s_j t_j + ... = s_j^2 t_j$. Since the mapping $t \mapsto D_1(s_j^{2}t)$ is unbounded, we have
 $\|D_1(s_j^2 t_j)\| \geq j + m \|\Lambda_{s_j}\|$, where $m$ is the bound of the bilinear mapping $(a, x)\mapsto xa:\mathcal{A} \times \mathcal{M} \rightarrow \mathcal{M}$. Now we have the following expressions:
\begin{align*}
\|G(s_j) H(c)\| & = \|D_1(s_j c) - E(s_j)F(c)\| \\ & \geq \|D_1(s_j c)\| - \|E(s_j)F(c)\| \\ & = \|D_1(s_j^2 t_j)\| - \|E(s_j)F(c)\| \\ & = \|D_1(s_j^2 t_j)\| - \|\Lambda_{s_j}(c)\| \\ & \geq
j + m \|\Lambda_{s_j}\| - m \|\Lambda_{s_j}\| \\ & = j.
\end{align*}
Since $\|s_j\| \leq 1$ and the mapping $t \mapsto G(t) H(c):\mathcal{A} \rightarrow \mathcal{M}$ is continuous, the non-equality $\|G(s_j) H(c)\| \geq j$ is a contradiction. This contradiction proves our claim that $\frac{\mathcal{A}}{\mathcal{I}}$ is finite-dimensional. \\
\\
\textbf{Step 5}: $D_1$ is continuous. \\ \\
Since the algebra $\frac{\mathcal{A}}{\mathcal{I}}$ is finite-dimensional, we can consider the elements $a_1, a_2, ..., a_r$ of $\mathcal{A}$ such that $\pi(a_1), \pi(a_2), ..., \pi(a_r)$ forms a basis for the algebra $\frac{\mathcal{A}}{\mathcal{I}}$. Suppose that $\lambda_1, \lambda_2, ..., \lambda_r$ are linear functionals on $\frac{\mathcal{A}}{\mathcal{I}}$ such that \[ \lambda_j(\pi(a_k)) =
\left\lbrace
  \begin{array}{c l}
     1 & \ \text{$j = k$ }\\
 0 & \text{ $ j \neq k$}
  \end{array}
\right. \]
As an easy exercise in functional analysis, we know that every $\lambda_j$ is continuous for $1 \leq j \leq r$. Since $\{\pi(a_1), \pi(a_2), ..., \pi(a_r)\}$ is a basis for the algebra $\frac{\mathcal{A}}{\mathcal{I}}$, for every element $a \in \mathcal{A}$ we have $\pi(a) = \sum_{j = 1}^{r}c_j \pi(a_j)$, where $c_j \in \mathbb{C}$. Hence, $\lambda_j(\pi(a)) = c_1 \lambda_j(\pi(a_1)) + c_2 \lambda_j(\pi(a_2)) + ... + c_j \lambda_j(\pi(a_j)) + ... + c_r \lambda_j(\pi(a_r)) = c_j$. Considering the continuous linear functionals $\varphi_j = \lambda_j o \pi$ ($1 \leq j \leq r$), we have
\begin{align*}
\pi(a) & = \sum_{j = 1}^{r}c_j \pi(a_j) \\ & = \sum_{j = 1}^{r}\lambda_j(\pi(a)) \pi(a_j) \\ &= \sum_{j = 1}^{r}\varphi_j(a) \pi(a_j).
\end{align*}

Consequently, $a - \sum_{j = 1}^{r}\varphi_j(a) a_j \in \mathcal{I}$. Now we define $\Phi:\mathcal{A} \rightarrow \mathcal{I}$ by $\Phi(a) = a - \sum_{j = 1}^{r}\varphi_j(a) a_j$. Obviously, the linear mapping $\Phi$ is continuous and so $D_1|_{\mathcal{I}} o \Phi:\mathcal{A} \rightarrow \mathcal{M}$ defined by $(D_1|_{\mathcal{I}} o \Phi)(a) = D_1(a - \sum_{j = 1}^{r}\varphi_j(a) a_j) = D_1(a) - \sum_{j = 1}^{r}\varphi_j(a) D_1(a_j)$ is continuous. The continuity of the mapping $D_1|_{\mathcal{I}} o \Phi$ along with the continuity of $\varphi_1, \varphi_2, ..., \varphi_r$ imply that the linear mapping \\ $a \mapsto D_1(a) - \sum_{j = 1}^{r}\varphi_j(a) D_1(a_j) + \sum_{j = 1}^{r}\varphi_j(a) D_1(a_j) = D_1(a):\mathcal{A} \rightarrow \mathcal{M}$ is continuous. Using a similar argument, we can prove that $D_1$ is continuous whenever $E(c) (F(ab) - F(a)F(b)) = 0$ for all $a, b, c \in \mathcal{A}$ and we leave it to the interested reader. Thereby, our proof is complete.
\end{proof}

In the following, there are some immediate consequences of the above theorem.

\begin{corollary} Let $\mathcal{A}$ be a $C^{\ast}$-algebra, let $\mathcal{B}$ be a Banach algebra and let $\mathcal{M}$ be a Banach $\mathcal{B}$-bimodule. Suppose that $E, H:\mathcal{A} \rightarrow \mathcal{M}$ are linear mappings and $F, G:\mathcal{A} \rightarrow \mathcal{B}$ are continuous at zero such that at least one of them is a homomorphism. Then every $(E, F, G, H)$-derivation from $\mathcal{A}$ into $\mathcal{M}$ is continuous.
\end{corollary}



If $\mathcal{Y}$ and $\mathcal{Z}$ are Banach spaces and $T: \mathcal{Y} \rightarrow \mathcal{Z}$ is a linear mapping, then the set $$S(T) = \left\{z \in \mathcal{Z} \ : \ \exists \ \{y_n\} \subseteq \mathcal{Y} \ such \ that \ y_n \rightarrow 0, T(y_n) \rightarrow z\right\}$$ is called the separating space of $T$. By the closed graph Theorem, $T$ is
continuous if and only if $S(T) = \{0\}$. For additional information about separating spaces, the
reader is referred to \cite{D}.

\begin{corollary} \label{5} Let $\mathcal{A}$ be a $C^{\ast}$-algebra, let $\mathcal{B}$ be a Banach algebra, and let $\mathcal{M}$ be a Banach $\mathcal{B}$-bimodule. Suppose that $E, H: \mathcal{A} \rightarrow \mathcal{M}$ are linear mappings and $F, G: \mathcal{A} \rightarrow \mathcal{B}$ are continuous mappings at zero such that $$\{m_0 \in \mathcal{M} \ : \ G(\mathcal{A}) m_0 = \{0\} \} = \{0\} = \{m_1 \in \mathcal{M} \ : \ m_1 F(\mathcal{A})= \{0\}\}.$$ If $D_1:\mathcal{A} \rightarrow \mathcal{M}$ is an $(E, F, G, H)$-derivation such that $(G(ab) - G(a)G(b))H(c) = 0$ or $E(c) (F(ab) - F(a)F(b)) = 0$ for all $a, b, c \in \mathcal{A}$, then $D_1$, $E$ and $H$ are continuous.
\end{corollary}

\begin{proof} According to Theorem \ref{4}, $D_1$ is continuous. Let $m_0 \in S(H) \subseteq \mathcal{M}$. Then there exists a sequence $\{b_n\} \subseteq \mathcal{A}$ such that $\lim_{n \rightarrow \infty}b_n = 0$ and $\lim_{n \rightarrow \infty}H(b_n) = m_0$. For arbitrary $a \in \mathcal{A}$, we have
\begin{align*}
0 = \lim_{n \rightarrow \infty}D_1(a b_n) = \lim_{n \rightarrow \infty}(E(a) F(b_n) + G(a) H(b_n)) = G(a)m_0,
\end{align*}
which means that $G(\mathcal{A}) m_0 = \{0\}$. By hypothesis, $m_0 = 0$ and this implies that $H$ is continuous. Similarly, we can show that $E$ is a continuous linear mapping, as desired.
\end{proof}

\begin{corollary} \label{*} Let $\mathcal{A}$ be a $C^{\ast}$-algebra, let $E, H: \mathcal{A} \rightarrow \mathcal{A}$ be linear mappings and let $F, G: \mathcal{A} \rightarrow \mathcal{A}$ be surjective and continuous mappings at zero. If $D_1:\mathcal{A} \rightarrow \mathcal{A}$ is an $(E, F, G, H)$-derivation such that $(G(ab) - G(a)G(b))H(c) = 0$ or $E(c) (F(ab) - F(a)F(b)) = 0$ for all $a, b, c \in \mathcal{A}$, then $D_1$, $E$ and $H$ are continuous.
\end{corollary}
\begin{proof} It is a well-known fact that every $C^{\ast}$-algebra is semiprime. Moreover, one can easily show that the only element $a_0$ of $\mathcal{A}$ which satisfies $\mathcal{A} a_0 = \{0\}$ or $a_0 \mathcal{A} = \{0\}$ is zero. Thereby, the previous corollary completes the proof.
\end{proof}



\begin{definition} Let $\sigma, \tau: \mathcal{A} \rightarrow \mathcal{B}$ be two mappings. A linear mapping $D:\mathcal{A} \rightarrow \mathcal{M}$ is called a generalized $(\sigma, \tau)$-derivation associated with a linear mapping $d:\mathcal{A} \rightarrow \mathcal{M}$ if $D(ab) = D(a)\sigma(b) + \tau(a) d(b)$ for all $a, b \in \mathcal{A}$.
\end{definition}
In the results presented on the continuity of generalized $(\sigma, \tau)$-derivations, it is assumed that the linear mapping $d$ is a $(\sigma, \tau)$-derivation. For example, see, \cite{Hu, J}. We remove this restriction in the next result.
\begin{corollary} Let $\mathcal{A}$ be a $C^{\ast}$-algebra, let $\mathcal{B}$ be a Banach algebra, let $\mathcal{M}$ be a Banach $\mathcal{B}$-bimodule and let $\sigma, \tau:\mathcal{A} \rightarrow \mathcal{B}$ be mappings. Let $D: \mathcal{A} \rightarrow \mathcal{M}$ be a generalized $(\sigma, \tau)$-derivation associated with a linear mapping $d:\mathcal{A} \rightarrow \mathcal{M}$ such that $\sigma, \tau:\mathcal{A} \rightarrow \mathcal{B}$ are continuous mappings at zero. Moreover, assume that either $(\tau(ab) - \tau(a)\tau(b))d(c) = 0$ or $D(c) (\sigma(ab) - \sigma(a)\sigma(b)) = 0$ for all $a, b, c \in \mathcal{A}$. Then $D$ is continuous.
\end{corollary}
\begin{proof} This is a direct consequence of Theorem \ref{4}.
\end{proof}




\bibliographystyle{amsplain}

\vskip 0.5 true cm

{\tiny (Amin Hosseini) Kashmar Higher Education Institute, Kashmar, Iran}
	
	{\tiny\textit{E-mail address:} hosseini.amin82@gmail.com}

	

\vskip 0.3 true cm 	
\end{document}